\documentclass[12pt]{article}
\usepackage{amsthm, amssymb, amsmath, latexsym}
\usepackage[utf8]{inputenc}
\usepackage[english]{babel}
\usepackage{longtable}

\begin{document}
\pagestyle{myheadings}
\thispagestyle{empty}
\setcounter{page}{1}

\newtheorem{definition}{Definition}
\newtheorem{proposition}{Proposition}
\newtheorem{theorem}{Theorem}
\newtheorem{lemma}{Lemma}
\newtheorem{corollary}{Corollary}
\newtheorem{remark}{Remark}
\theoremstyle{plain}
\mathsurround 2pt

\gdef\Aut{\mathop{\rm Aut}\nolimits}
\gdef\End{\mathop{\rm End}\nolimits}
\gdef\Ker{\mathop{\rm Ker}\nolimits}
\gdef\Im{\mathop{\rm Im}\nolimits}
\gdef\Inn{\mathop{\rm Inn}\nolimits}
\gdef\exp{\mathop{\rm exp}\nolimits}

\begin{center}
\textbf{\Large Semidistributive nearrings with identity}
\end{center}

\begin{center}
{\small Iryna Raievska$^{2,1}$, Maryna Raievska$^{2,1}$, Yaroslav Sysak$^1$\\
${}^1$Institute of Mathematics of NAS of Ukraine, Ukraine\\
${}^2$University of Warsaw, Poland\\
raeirina@imath.kiev.ua, raemarina@imath.kiev.ua, sysak@imath.kiev.ua}
\end{center}

\begin{abstract}
It is proved that the additive group of every semidistributive nearring R with an identity is abelian and if R has no elements of order $2$, then the nearring R actually is an associative ring
\end{abstract}

\section{Introduction}

Nearrings are a generalization of associative rings in the sense that with the respect to addition they need not be commutative and only one (left or right) distributive law is assumed. In this paper the concept ``nearring'' means a left distributive nearring. Such a nearring $R$ is called a semidistributive if $(r+s+r)t=rt+st+rt$ for all $r,~s,~t\in R$. Clearly every associative ring is a semidistributive nearring, but not conversely. It is proved that the additive group of every semidistributive nearring $R$ with an identity is abelian and if $R$ has no elements of order $2$, then the nearring $R$ actually is an associative ring. The reader is referred to the books by Meldrum~\cite{Meldrum_1985} or Pilz~\cite{Pilz_1977} for terminology, definitions and basic facts concerning nearrings.

\section{Preliminaries}

We recall first some basic definitions of the theory of nearrings.

\begin{definition}
A set $R$ with two binary operations $``+"$ and $``\cdot"$ is called a (left) nearring if the following statements hold:
\begin{itemize}
  \item[(1)] $(R,+)$ is a (not necessarily abelian) group with neutral element $0$;
  \item[(2)] $(R,\cdot)$ is a semigroup;
  \item[(3)] $x\cdot (y+z)=x\cdot y+x\cdot z$ for all $x$, $y$, $z\in R$.
\end{itemize}
\end{definition}

If $R$ is a nearring, then the group $R^+=(R,+)$ is called the additive group of $R$. Furthermore, if $M$ is a subgroup of $R^+$, then it follows from statement~(3) that for each element $x\in R$ the set $xM=x\cdot y| y\in M\}$ is a subgroup of $R^+$ and, in particular, $x \cdot 0 = 0$. If in addition $0\cdot x=0$, then the nearring R is called zero-symmetric and if the semigroup $(R, \cdot)$ is a monoid, i.~e. it has an identity element $i$, then $R$ is a nearring with identity $i$. In the latter case the group $R^*$ of all invertible elements of the monoid $(R, \cdot)$ is called the multiplicative group of $R$. A subgroup $M$ of $R^+$ is called $R^*$-invariant if $rM \leq M$ for each $r\in R^*$ and $M$ is an $(R,R)$-subgroup, if $xMy\subseteq M$ for arbitrary elements $x,~y\in R$.

It is well-known that the set of all endomorphisms of a group $G$ forms a monoid with respect to composition of endomorphisms. This monoid will be denoted by $\End(G)$ and called the endomorphism monoid of $G$. Clearly the automorphism group $\Aut(G)$ is a submonoid of $\End(G)$ consisting of the set of all invertible elements of $\End(G)$. The following lemma establishes a connection between some semigroups of endomorphisms of the additive group of a nearring and its multiplicative subsemigroups.

\begin{lemma}\label{Lemma_1}
Let $R$ be a nearring with identity $i$ and $S$ a semigroup of $(R, \cdot)$. Then the endomorphism monoid $\End(R^+)$ contains a subsemigroup $T$ isomorphic to $S$ and satisfying the condition
$$i^T=\{i^t| t\in T\}=S.$$
In particular, if $S=(R, \cdot)$, then the intersection $A=\Aut(R+)\cap T$ is a subgroup of $\Aut(R+)$ isomorphic to the multiplicative group $R^*$ and
$$i^A=\{i^a| a\in A\}=R^*.$$
\end{lemma}

The subgroup $A$ of the automorphism group $\Aut(R^+)$ defined in Lemma~\ref{Lemma_1} is called \emph{the subgroup} of $\Aut(R^+)$ \emph{associated with the group} $R^*$.

\begin{definition}
A (left) nearring $R$ is called semidistributive if so is the multiplication from the right in respect to its addition. In other words, for any elements $r,~s,~t\in R$ the equality $(r+s+r)t=rt+st+rt$ holds.
\end{definition}

The following assertion is well-known (see, for instance, \cite{ClMal_66}, Theorem 3).

\begin{lemma}
The exponent of the additive group of a finite nearring $R$ with identity $i$ is equal to the additive order of $i$ which coincides with the additive order of every invertible element of $R$.
\end{lemma}

\section{Additive groups of semidistributive nearrings}

It is obvious that every distributive nearring is semidistributive, but not conversely. For example, the nearring $Map(G)$ of all functions on the group $G$ of order $2$ is semidistributive and not distributive.

\begin{lemma}\label{Lemma_2}
Let $R$ be a semidistributive nearring, $n$ a positive integer and let $r,~s$ be any elements of $R$. Then the following statements hold:
\begin{itemize}
  \item[(1)] $(-r)\cdot s=-(r\cdot s)$ and so either $0\cdot s=0$ or $0\cdot s$ is an element of order $2$ in $R^+$. In particular, if $r$ is an element of odd order in $R^+$, then $0\cdot r=0$;
  \item[(2)] $(r\cdot n)s=(rs)\cdot n$ for odd $n$ and $(r\cdot n)s=(rs)\cdot n+0\cdot s$ for even $n$;
  \item[(3)] if $r$ and $s$ are elements of coprime orders, then $r\cdot s=0$.
\end{itemize}
\end{lemma}

\emph{Proof.} Since $r\cdot s=(r+(-r)+r)\cdot s=r\cdot s+(-r)\cdot s+r\cdot s$, it follows that $0=(-r)\cdot s+r\cdot s$ whence $(-r)\cdot s=-(r\cdot s)$. In particular, if $r\cdot m=0$ for some odd $m$, then $1=2n-m$ for $n=\frac{m+1}{2}$ and so $0\cdot r=(0\cdot r)\cdot {2n}-(0\cdot r)\cdot m=0$, as claimed.

For $n=1$ statement~(2) is obvious. If $n>1$, then $(r\cdot n)s=(r+r\cdot (n-2)+r)s=rs+(r\cdot (n-2))s+rs=(rs)\cdot 2+(r\cdot (n-2))s$. By induction on $n$, we have $(r\cdot (n-2))s=(rs)\cdot (n-2)$ for odd $n$ and $(r\cdot (n-2))s=(rs)\cdot (n-2)+0\cdot s$ for even $n$. Therefore $(r\cdot n)s=(rs)\cdot 2+(rs)\cdot (n-2)=(rs)\cdot n$ for odd $n$ and $(r\cdot n)s=(rs)\cdot 2+(rs)\cdot (n-2)=(rs)\cdot n+0\cdot s$ for even $n$, proving~(2).

Finally, let $m$ and $n$ be the orders of $r$ and $s$, respectively. If $(m,n)=1$, then one of these numbers, say $m$, is odd and so $(rs)\cdot m=(r\cdot m)s=0\cdot s=0$ by statement~(2). On the other hand, $(rs)\cdot n=r(s\cdot n)=r\cdot 0=0$. Therefore $rs=0$, as desired.

It is well-known that the additive group of any distributive nearring with identity is abelian. The following assertion extends this fact to semidistributive nearrings.

\begin{lemma}\label{Lemma_3}
The additive group of every semidistributive nearring $R$ with an identity is abelian.
\end{lemma}

\emph{Proof.} Indeed, if $i$ is an identity of $R$, then $(-i)s=-s$ for each $s\in R$ by statement~(1) of Lemma~\ref{Lemma_2} and therefore $(-i)(r+s)=-(r+s)=-s-r$. On the other hand, $(-i)(r+s)=(-i)r+(-i)s=-r-s$ by the left distributivity. Hence $-s-r=-r-s$ and thus $r+s=s+r$, as claimed.

On the other hand, the additive groups of semidistributive nearrings without unity do not have to abelian, as the following example shows.

{\bf Example.} Let $N$ be a nearring whose additive group $N^+$ is isomorphic to a symmetric group of degree 3. Then $N^+=\langle a\rangle +\langle b\rangle$ with $a\cdot 3=b\cdot 2=0$ and $b+a+b=a\cdot 2$. In particular, each element of $N^+$ can be uniquely written in the form $a\cdot m+b\cdot n$ with $0\leq m\leq 2$ and $0\leq n\leq1$. We define a multiplication $*$ on $N^+$ by the table

\

\begin{tabular}{ccccccc}
    $*$          & 0 & $a$ & $a\cdot 2$ & $b$          & $a+b$        & $a\cdot 2+b$ \\
    0            & 0 & 0   & 0          & 0            & 0            & 0 \\
    $a$          & 0 & 0   & 0          & 0            & 0            & 0 \\
    $a\cdot 2$   & 0 & 0   & 0          & 0            & 0            & 0 \\
    $b$          & 0 & 0   & 0          & $b$          & $b$          & $b$ \\
    $a+b$        & 0 & 0   & 0          & $a+b$        & $a+b$        & $a+b$ \\
    $a\cdot 2+b$ & 0 & 0   & 0          & $a\cdot 2+b$ & $a\cdot 2+b$ & $a\cdot 2+b$ \\
  \end{tabular}

\

Looking at this table, it is easy to see that $(N,+,*)$ is a nearring that is semidistributive and zero-symmetric.

We note in passing that the package SONATA \cite{SONATA} of the computer algebra system GAP \cite{GAP} contains a library of all nearrings of order at most 15. In particular, there exist 39 non-isomorphic nearrings whose additive groups are isomorphic to the symmetric group of degree 3, among which only 4 are semi-distributive, including 2 distributive.

\begin{lemma}\label{Lemma_4}
Let $R$ be a semidistributive nearring with an identity. If the additive group $R^+$ is torsion-free, then $R$ is a ring.
\end{lemma}

\emph{Proof.} Clearly it suffices to prove that in $R$ the right distributive law holds. Indeed, since the additive group $R^+$ is abelian by Lemma~\ref{Lemma_3}, for any elements $r,~s,~t$ of $R$ we have $(r+s)t=(r+s+(-(r+s))+s+r)t=rt+(s+(-(r+s))+s)t+rt=rt2+st2+(-(r+s))t$ and $(-(r+s))t=-(r+s)t$ by statement~(1) of Lemma~\ref{Lemma_2}. Therefore, $((r+s)t-rt-st)2=0$ whence $(r+s)t=rt+st$, as claimed.

\section{Distributive elements in semidistributive\\
nearrings}

Recall that an element $t$ of a nearring $R$ is called distributive in $R$ if $(r+s)t=rt+st$ for any elements $r,~s$ of $R$.

\begin{lemma}\label{Lemma_5}
Let $R$ be a semidistributive nearring with an identity. Then the elements of odd orders of the additive group of $R$ are distributive in $R$. In particular, each semidistributive nearring of odd order is a ring.
\end{lemma}

\emph{Proof.} Let $t$ be an element of odd order $n$ in $R^+$ and $m=\frac{n+1}{2}$. Then $u=tm$ is an element of order $n$ and $t=u+u$. Since the nearring $R$ is semidistributive, for any elements $r,s\in R$ we have $(r+s)u=(r+(s-r)+r)u=ru+(s-r)u+ru$ and similarly $(s+r)u=(s+(r-s)+s)u=su+(r-s)u+su$. Furthermore, as $r+s=s+r$ by Lemma~\ref{Lemma_3} and $(r-s)u=(-(s-r))u=-(s-r)u$ by statement~(1) of Lemma~\ref{Lemma_2}, it follows that $(r+s)u=ru+(s-r)u+ru=su-(s-r)u+su$. Therefore $(r+s)t=(r+s)(u+u)=r(u+u)+s(u+u)=rt+st$ and so the element $t$ is distributive. In particular, if $R$ is of odd order, then $R$ is a distributive nearring with identity and so a ring.

We recall that the concepts of a \emph{subnearring} and a \emph{nearring homomorphism} are defined by the same way as for rings. In particular, if $\lambda$ is a nearring homomorphism of $R$, then its kernel $\Ker \lambda$ is a subnearring of R whose additive subgroup is normal in $R^+$. A subnearring $I$ of $R$ is an ideal of $R$ if $I=\Ker \lambda$ for some $\lambda$. It can simply be verified that $I$ is an ideal of $R$ if and only if $I^+$ is a normal subgroup of $R^+$ and for any elements $r,~s\in R$ and $a\in I$ the inclusions $ra\in I$ and $(r+a)s-rs\in I$ hold.

\begin{lemma}\label{Lemma_6}
Let $R$ be a semidistributive nearring with identity. If $p$ is a prime number and $P$ is the $p$-component of the additive group $R^+$, then $P$ is an ideal of $R$.
\end{lemma}

\emph{Proof.} Since the additive group $R^+$ is abelian by Lemma~\ref{Lemma_3}, every $p$-element of $R$ belongs to $P$ and for each $a\in P$ there exists a positive integer $n$ such that $ap^n=0$. As $(ra)p^n=r(ap^n)=r0=0$ for every $r\in R$, this implies $ra\in P$. Thus, in order to prove that $P$ is an ideal of $R$, it remains to show that for each $s\in R$ the element $x=(r+a)s-rs$ belongs to $P$.

Assume first that $p$ is odd, so that $a$ is of odd order. Then there exists an element $b\in P$ of the same order such that $a=b+b$. Putting $t=s+s$, we have $x=(b+r+b)s-rs=bs+rs+bs-rs=bs+bs=bt$ and so $xp^n=(bt)p^n=(bp^n)t=0t=0$ by statement~(2) of Lemma~\ref{Lemma_3}. Therefore $x=0\in P$, as desired.

Now let $p=2$ and so $a2^{n}=0$. Then $x2^{n}=((r+a)s-rs)2^{n}=((r+a)s)2^{n}-(rs)2^{n}$. Further, applying statement~(2) of Lemma~\ref{Lemma_2} we deduce $((r+a)s)2^{n}+0s=((r+a)2^{n})s=(r2^{n}+a2^{n})s=(r2^{n})s=(rs)2^{n}+0s$. From this it follows that $((r + a)s)2^{n}=(rs)2^{n}$ and hence $x2^{n}=((r+a)s)2^{n}-(rs)2^{n}=0$. Therefore, $x\in P$ and this completes the proof.

\section{Modules over nearrings and semidistributivity}

Let $R$ be a nearring and $G$ an additive group. Then $G$ is said to be a (right) $R$-\emph{module} if for any elements $g\in G$ and $r\in R$ there exists a unique element $gr\in G$ such that $(gr)s=g(rs)$ and $g(r+s)=gr+gs$ for every $s\in R$. It is clear that the additive group of a nearring $R$ is an $R$-module over $R$ with respect to nearring multiplication from the right as the action. In addition, if $H$ is a subset of an $R$-module $G$, then the set $Ann_R(H)=\{x\in R| Hx=0\}$ is called the (right) annihilator of $H$ in $R$. An $R$-module $G$ is called \emph{faithful} if $Ann_R(G)=0$.

\begin{lemma}\label{Lemma_7}
If $R$ is a nearring and $G$ is an $R$-module, then the annihilator $AnnR(G)$ is an ideal of $R$.
\end{lemma}

\emph{Proof.} Indeed, let $g\in G$ and $x\in Ann_R(G)$. Since $gr\in G$ for each $r\in R$, it follows that $g(rx)=(gr)x=0$. Therefore $rx\in Ann_R(G)$ and hence $Ann_R(G)$ is a left ideal of $R$. To prove that $Ann_R(G)$ is an ideal of $R$, it is enough to show that for any $r,~s\in R$ the element $t=(r+x)s-rs$ belongs to $Ann_R(G)$. As $G$ is an $R$-module and $gx=0$, we have $g(r+x)=gr+gx=gr$ and so $gt=g((r+x)s)+g(-rs)=(g(r+x))s+g(-rs)=g(rs)+g(-rs)=g0=0$. This implies $t\in Ann_R(G)$ and thus $Ann_R(G)$ is an ideal of $R$, as claimed.

The following lemma is proved in \cite{Meldrum_1985}, Lemma~3.7.

\begin{lemma}\label{Lemma_8}
Let $R$ be a nearring and $G$ an $R$-module with the properties:
\begin{itemize}
  \item[(i)] $G$ is faithful as an $R$-module,
  \item[(ii)] $(G, +)$ is an abelian group,
  \item[(iii)] for each $r\in R$ the mapping $\hat{r} \colon G \rightarrow G$ given by $\hat{r} \colon g \rightarrow gr$ is an endomorphism of $G$.
\end{itemize}
Then R is a ring.
\end{lemma}

\emph{Proof.} Let $s,~t$ be elements of $R$ and $r=s+t$. Then $gr=g(s+t)=gs+gt =gt+gs=g(t+s)$ by property $(ii)$. Furthermore, from $(iii)$ it follows that $g(s+t)r =(gs+gt)^{\hat{r}}=(gs)r+(gt)r=g(sr+tr)$. Finally, since the $R$-module $G$ is faithful by $(i)$, this implies $(s+t)r=sr+tr$ and so $R$ is a distributive ring.

\begin{theorem}
A semidistributive nearring $R$ with identity is simple as a nearring if and only if $R$ is a simple associative ring.
\end{theorem}

\emph{Proof.} As noted above, the additive group $G=R^+$ can be viewed as an $R$-module over $R$ with respect to nearring multiplication from the right as the action. If the nearring $R$ is simple, it has no proper ideals and so $Ann(R(G))=0$ by Lemma~\ref{Lemma_7}. Therefore, $G$ is faithful as an $R$-module and $G=R^+$ is an abelian group by Lemma~\ref{Lemma_3}. Finally, for each $r\in R$ the mapping $\rho \colon G \rightarrow G$ given by $\rho \colon g\rightarrow gr$ is an endomorphism of $G$, because $\rho(g+h)=r(g+h)=rg+rh=\rho(g)+\rho(h)$ for every $g,~h\in G$. Applying now Lemma~\ref{Lemma_8}, we deduce that $R$ is a simple ring. The converse is obvious.

\

\footnotesize{{\bf Acknowledgement.} The first and second authors would like to thank to IIE-SRF for supporting of their fellowship at the University of Warsaw. The third author is grateful for the financial support within the program of support for priority research and technical (experimental) developments of the Section of Mathematics of the NAS of Ukraine for 2022-2023. Project “Innovative methods in the theory of differential equations, computational mathematics and mathematical modeling”, No. 7/1/241. (State Registration No. 0122U000670).}

\end{document}